\input amstex
\input amsppt.sty   
\pagewidth{11.8cm} 
\pageheight{17.5cm}
\magnification=\magstep1
\def\nmb#1#2{#2}         
\def\ign#1{}             

\redefine\o{\circ}
\define\X{\frak X}
\define\al{\alpha}

\define\de{\delta}
\define\ep{\varepsilon}
\define\ze{\zeta}
\define\et{\eta}

\define\ka{\kappa}
\define\la{\lambda}

\define\ph{\varphi}

\define\ps{\psi}

\define\Ga{\Gamma}

\define\Ph{\Phi}

\predefine\ii{\i}
\redefine\i{^{-1}}
\define\row#1#2#3{#1_{#2},\ldots,#1_{#3}}
\define\x{\times}
\define\M{\Cal M}
\redefine\D{\Cal D}
\define\Vol{\operatorname{vol}}
\define\Diff{\operatorname{Diff}}
\define\tr{\operatorname{tr}}
\define\dd#1{\frac{\partial}{\partial#1}}
\define\Exp{\operatorname{Exp}}
\define\arctg{\operatorname{arctg}}

\def\today{\ifcase\month\or
 January\or February\or March\or April\or May\or June\or
 July\or August\or September\or October\or November\or December\fi
 \space\number\day, \number\year}
 
\topmatter
\title The Riemannian Manifold of all Riemannian Metrics \endtitle
\author  Olga Gil-Medrano, 
Peter W. Michor
\endauthor
\address{Departamento de Geometr\'\ii a y Topolog\'\ii a,
Facultad de Matem\'aticas,
Universidad de Valencia,
46100 Burjasot,
Valencia, Spain.}\endaddress
\address{Institut f\"ur Mathematik, Universit\"at Wien,
Strudlhofgasse 4, A-1090 Wien, Austria.}
\endaddress
\thanks{The first author was partially supported the CICYT grant
n. PS87-0115-G03-01. This paper was prepared during a stay of the second
author in Valencia, by a grant given by Conseller\'\ii a de Cultura,
Educaci\'on y Ciencia, Generalidad Valenciana.}\endthanks
\endtopmatter

\document
 
\heading Introduction \endheading
 
 If $M$ is a (not necessarily compact)
smooth finite dimensional manifold, the space $\M=\M(M)$
of all Riemannian metrics on it can be endowed with
a structure of an infinite dimensional smooth manifold modeled
on the space $\D(S^2T^*M)$ of symmetric $\binom02$-tensor fields with
compact support, in the sense of \cite{Michor, 1980}.
The tangent bundle of $\M$ is $T\M =
C^\infty(S^2_+T^*M)\x\D(S^2T^*M)$ and a smooth Riemannian
metric can be defined by
$$G_g(h,k)= \int_M\tr(g\i hg\i k)\Vol(g).$$
In this paper we study the geometry of $(\M,G)$ by using
the ideas developed in \cite{Michor, 1980}.
 
With that differentiable structure on $\M$ it is possible
to use variational principles and so we start in
section 2 by computing  geodesics as the curves in
$\M$ minimizing the energy functional. From the geodesic
equation, the covariant derivative of the Levi-Civita
connection can be obtained, and that provides a direct
method for computing the curvature of the manifold.
 
Christoffel symbol and curvature turn out to be
pointwise in M and so, although the mappings involved in
the definition of the Ricci tensor and the scalar
curvature have no trace, in our case we can define the
concepts of "Ricci like curvature" and "scalar like
curvature".
 
The pointwise character mentioned above allows us in
section 3,  to solve explicitly the geodesic equation
and to obtain the domain of definition of the
exponential mapping. That domain turns out to be open
for the topology considered on $\M$ and the
exponential mapping is a diffeomorphism onto its image
which is also explicitly given. In the $L^2$-topology given by $G$
itself this domain is,
however, nowhere open. Moreover, we prove
that it is, in fact, a real analytic diffeomorphism,
using \cite{Kriegl-Michor, 1990}.
We think that this exponential mapping will be a very powerful tool
for further investigations of the stratification of
orbit space of $\M$ under the diffeomorphism group,
and also the stratification of principal connections modulo the gauge
group.
 
In section 4 Jacobi fields of an infinite dimensional
Riemannian manifold are defined as the infinitesimal
geodesic variations and we show that they must
satisfy the Jacobi Equation. For the manifold
$(\M,G)$ the existence of Jacobi fields, with any
initial conditions, is obtained from the results about
the exponential mapping in section 3. Uniqueness and the
fact that they are exactly the solutions of the
Jacobi Equation follows from its pointwise character. We
finally give the expresion of the Jacobi fields.
 
For fixed $x\in M$, there exists a family of homothetic
Riemannian metrics in the finite dimensional manifold
$S^2_+T_x^*M$ whose geodesics are given by the
evaluation of the geodesics of $(\M,G)$. The relationship
between the geometry of $(\M,G)$ and that of these
manifolds is explained in each case and it is used to
visualize the exponential mapping. Nevertheless, in this
paper, we have not made use of these manifolds to obtain
the results, every computation having been made directly
on the infinite dimensional manifold.
 
Metrics on $S^2_+T^*_xM$ for three dimensional manifolds $M$ which
are similar to ours but have different signatures were
considered by \cite{DeWitt, 1967}. He computed the curvature and the
geodesics and gave some ideas on how to use them to determine the
distance between two 3-geometries, but without considering explicitly
the infinite dimensional manifold of all Riemannian metrics on a given
manifold.
 
The topology of $\M(M)$, under the assumption that $M$ is
compact, orientable, without boundary, was studied by
\cite{Ebin, 1970} who treated $G$ in the context of Sobolev
completions of mapping spaces and computed the Levi-Civita
connection. In the same context and under the
same assumptions, the curvature and the geodesics have been
computed in \cite{Freed-Groisser, 1989}. 

The explicit formulas of the three papers just mentioned are the same 
as in this paper.
 
We want to thank A. Montesinos Amilibia for
producing the computer image of figure 1.
 
\heading \nmb0{1}. The general setup \endheading
 
\subheading{\nmb.{1.1}. The space of Riemannian metrics}
Let $M$ be a smooth second countable
finite dimensional manifold. Let $S^2T^*M$ denote the vector bundle of all
symmetric $\binom02$-tensors on $M$ and let $S^2_+T^*M$ be the
open subset of all the positive definite ones. Then the space
$\M(M)=\M$ of all Riemannian metrics is the space of sections
$C^\infty(S^2_+T^*M)$ of this fiber bundle. It is open in the
space of sections $C^\infty(S^2T^*M)$ in the Whitney
$C^\infty$-topology, in which the latter space is, however, not
a topological vector space, since $\frac1n h$ converges to 0 if
and only if $h$ has compact support. So the space $\D=\D(S^2T^*M)$
of sections with compact support is the largest topological
vector space contained in the topological group
$(C^\infty(S^2T^*M),+)$, and the trace of the Whitney
$C^\infty$-topology on it coincides with the inductive limit topology
$$\D(S^2T^*M)=\varinjlim_KC^\infty_K(S^2T^*M),$$
where $C^\infty_K(S^2T^*M)$ is the space of all sections with
support contained in $K$ and where $K$ runs through all compact
subsets of $M$.
 
So we declare the path components of $C^\infty(M,S^2_+T^*M)$ for
the Whitney $C^\infty$-topology also to be open. We get a
topology which is finer than the Whitney topology, where each
connected component is homeomorphic to an open subset in
$\D=\D(S^2T^*M)$. So $\M=C^\infty(S^2_+T^*M)$ is a smooth manifold
modeled on nuclear (LF)-spaces, and the tangent bundle is given
by $T\M=\M\x\D$.
 
\subheading{\nmb.{1.2}. Remarks}
The main reference for the infinite dimensional manifold
structures is \cite{Michor, 1980}. But the
differential calculus used there is not completely up to date,
the reader should consult \cite{Fr\"olicher-Kriegl, 1988}, whose
calculus is more natural and much easier to apply.  There
a mapping between locally convex spaces is smooth if and only if
it maps smooth curves to smooth curves. See also
\cite{Kriegl-Michor, 1990} for a setting for real analytic
mappings along the same lines and applications to manifolds of
mappings.
 
As a final remark let us add that the differential structure on
the space $\M$ of Riemannian metrics is not completely
satisfying, if $M$ is not compact. In fact $C^\infty(S^2T^*M)$ is
a topological vector space with the compact $C^\infty$-topology,
but the space $\M=C^\infty(S^2_+T^*M)$ of Riemannian metrics is
not open in it. Nevertheless, we will see later that the
exponential mapping for the natural Riemannian metric on $\M$ is
defined also for some tangent vectors which are not in $\D$.
This is an indication that the most natural setting for
manifolds of mappings is based on the compact
$C^\infty$-topology, but that one loses existence of charts. In
\cite{Michor, 1985} a setting for infinite dimensional manifolds
is presented which is based on an axiomatic structure of smooth
curves instead of charts.
 
\subheading{\nmb.{1.3}. The metric} The tangent bundle of the
space $\M=C^\infty(S^2_+T^*M)$ of Riemannian metrics is
$T\M=\M\x\D=C^\infty(S^2_+T^*M)\x\D(S^2T^*M)$.
We identify the vector bundle $S^2T^*M$ with the subbundle
$$\{\ell\in L(TM,T^*M): \ell^t=\ell\}$$
of $L(TM,T^*M)$, where the transposed is given by the composition
$$\ell^t: TM @>i>> T^{**}M @>\ell^*>> T^*M.$$
Then the fiberwise inner product on $S^2T^*M$ induced by $g\in\M$
is given by the expression
$\langle h,k\rangle_g:= \tr(g\i hg\i k)$, so a smooth Riemannian
metric on $\M$ is given by
$$G_g(h,k)=\int_M\tr(g\i hg\i k)\Vol(g),$$
where $\Vol(g)$ is the positive density defined by the local formula
$\Vol(g)=\sqrt{\det g}\,dx$.
We call this the {\it canonical Riemannian metric} on $\M$,
since it is invariant under the action of the diffeomorphism
group $\Diff(M)$ on the space $\M$ of metrics. The integral is
defined since $h$, $k$ have compact support. The metric is
positive definite, $G_g(h,h)\geq0$ and $G_g(h,h)=0$ only if
$h=0$. So $G_g$ defines a linear injective mapping
from the tangent space $T_g\M=\D(S^2T^*M)$ into its dual
$\D(S^2T^*M)'$, the space of distributional densities with values
in the dual bundle $S^2TM$. This linear mapping is, however, never
surjective, so $G$ is only a {\it weak} Riemannian metric. The
tangent space $T_g\M=\D(S^2T^*M)$ is a pre-Hilbert space, whose
completion is a Sobolev space of order 0, depending on $g$ if
$M$ is not compact.
 
\subheading{\nmb.{1.4}. Remark} Since $G$ is only a weak
Riemannian metric, all objects which are only implicitly given,
a priori lie in the Sobolev completions of the relevant spaces.
In particular this applies to the formula
$$\align
2G(\xi,\nabla_\et\ze)=&\xi G(\et,\ze) + \et G(\ze,\xi) - \ze G(\xi,\et)\\
&+ G([\xi,\et],\ze) +G([\et,\ze],\xi) - G([\ze,\xi],\et),
\endalign$$
which a priori gives only uniqueness but not existence of the
Levi Civita covariant derivative.
 
\heading \nmb0{2}. Geodesics, Levi Civita connection,\\ and curvature
\endheading
 
\subheading{\nmb.{2.1} The covariant derivative} Since we will
need later the covariant derivative of vector fields along a
geodesic for the derivation of the Jacobi equation, we present
here a careful description of the notion of the covariant
derivative, which is valid in infinite dimensions.
Here $\M$ might be any infinite dimensional manifold, modeled on
locally convex spaces. If we are given a horizontal bundle,
complementary to the vertical one, in $T^2\M$, with the usual
properties of a linear connection, then the projection from
$T^2\M$ to the vertical bundle $V(T\M)$ along the horizontal
bundle, followed by the vertical projection $V(T\M)\to T\M$,
defines the {\it connector} $K:T^2\M \to T\M$, which has the
following properties:
\roster
\item It is a left inverse to the vertical lift mapping with any
        foot point.
\item It is linear for both vector bundle structures on $T^2M$.
\item The connection is symmetric (torsionfree) if and only if
        $K\o \ka = K$, where $\ka$ is the canonical flip mapping
        on the second tangent bundle.
\endroster
If a connector $K$ is given, the covariant derivative is defined
as follows: Let $f:\Cal N\to \M$ be a smooth mapping, let
$s:\Cal N\to T\M$ be a vector field along $f$ and let
$X_x\in T_x\Cal N$.
Then
$$\nabla_{X_x}s := (K\o Ts)(X_x).$$
In a chart the Christoffel symbol is related to the connector by
$$K(g,h;k,\ell) = (g,\ell-\Ga_g(h,k))$$.
 
We want to state one property, which is usually stated rather
clumsily in the literature: If $f_1:\Cal P\to \Cal N$ is another smooth
mapping and $Y_y\in T_y\Cal P$, then we have
$\nabla_{Y_y}(s\o f_1) = \nabla_{T_y(f_1)Y_y}s$. Equivalently,
if vector fields $Y \in \X(\Cal P)$ and $X \in \X(\Cal N)$ are
$f_1$-related, then $\nabla_Y(s\o f_1)=(\nabla_Xs)\o f_1$.
 
If $V(t)$ is a vector field along a smooth curve
$g(t)$, we have $\nabla_{\partial_t}V = \dd t V - \Ga_g(g_t,V)$,
in local coordinates.
 
If $c: \Bbb R^2 \to \M$ is a smooth mapping for a symmetric
connector $K$ we have
$$\multline
\nabla_{\partial_t}\dd sc(t,s) = K\o T(Tc\o \partial_s)\o \partial_t
= K\o T^2c\o T(\partial_s)\o \partial_t  \\
= K\o\ka\o T^2c\o T(\partial_s)\o \partial_t
= K\o T^2c\o\ka\o T(\partial_s)\o \partial_t \\
= K\o T^2c\o T(\partial_t)\o \partial_s
= \nabla_{\partial_s}\dd tc(t,s),
\endmultline$$
which will be used for Jacobi fields.
 
\subheading{\nmb.{2.2}} Let $t\mapsto g(t)$ be a smooth curve in
$\M$: so $g:\Bbb R\x M\to S^2_+T^*M$ is smooth and by the choice of
the topology on $\M$ made in \nmb!{1.1} the curve
$g(t)$ varies only in a compact subset of $M$, locally
in $t$, by \cite{Michor, 1980, 4.4.4, 4.11, and 11.9}.
Then its energy is given by
$$\align E_a^b(g):&= \tfrac12\int_a^bG_g(g_t,g_t)dt\\
&=\tfrac12\int_a^b\int_M\tr(g\i g_t g\i g_t)\Vol(g)\,dt,
\endalign$$
where $g_t=\frac{\partial}{\partial t}g(t)$.
 
Now we consider a variation of this curve, so we assume now that
$(t,s)\mapsto g(t,s)$ is smooth in all variables and locally in
$(t,s)$ it only varies within a compact subset in $M$ --- this is
again the
effect of the topology chosen in \nmb!{1.1}. Note that $g(t,0)$ is
the old $g(t)$ above.
 
\proclaim{\nmb.{2.3}. Lemma} In the setting of \nmb!{2.2} we
have the first variation formula
$$\multline \frac{\partial}{\partial s}|_0E_a^b(g(\quad,s))
= G_g(g_t,g_s)|_{t=a}^{t=b} +\\
+ \int_a^bG_g\bigl(
-g_{tt}+ g_t g\i g_t + \tfrac14\tr(g\i g_t g\i g_t)g -\tfrac12\tr(g\i g_t)g_t
        ,g_s\bigr)\,dt.
\endmultline$$
\endproclaim
\demo{Proof} We have
$$\frac{\partial}{\partial s}|_0E_a^b(g(\quad,s))
=\frac{\partial}{\partial s}|_0
\frac12\int_a^b\int_M\tr(g\i g_t g\i g_t)\Vol(g)dt.$$
We may interchange $\frac{\partial}{\partial s}|_0$ with the
first integral since this is finite dimensional analysis, and we
may interchange it with the second one, since $\int_M$ is a
continuous linear functional on the space of all smooth
densities with compact support on $M$, by the chain rule.
Then we use that $\tr_*$ is linear and continuous,
$d(\Vol)(g)h= \frac12\tr(g\i h)\Vol(g)$, and that
$d((\quad)\i)_*(g)h= -g\i h g\i$ and partial integration.
\qed\enddemo
 
\subheading{\nmb.{2.4} The geodesic equation}
By lemma \nmb!{2.3} the curve $t\mapsto g(t)$ is a geodesic if
and only if we have
$$\align g_{tt} &= g_t g\i g_t + \tfrac14 \tr(g\i g_t g\i g_t)\,g
        - \tfrac12 \tr(g\i g_t)\,g_t\\
&= \Ga_g(g_t,g_t),
\endalign$$
where the {\it Christoffel symbol} $\Ga:\M\x\D\x\D\to \D$ is
given by symmetrisation
$$\multline \Ga_g(h,k) = \tfrac12 h g\i k + \tfrac12 k g\i h + \\
+ \tfrac14 \tr(g\i h g\i k)\,g
- \tfrac14 \tr(g\i h)\,k  - \tfrac14 \tr(g\i k)\,h.
\endmultline$$
The sign of $\Ga$ is chosen in such a way that the horizontal
subspace of $T^2\M$ is parameterized by $(x,y;z,\Ga_x(y,z))$.
If instead of the obvious framing we use
$T\M=\M\x\D \ni (g,h)\mapsto (g,g\i h)=:(g,H)\in \{g\}\x \D(L_{sym,g}(TM,TM))
\subset \M\x\D(L(TM,TM))$, the Christoffel symbol looks like
$$\bar\Ga_g(H,K) =
\tfrac12(HK+KH) + \tfrac14 \tr(HK)Id - \tfrac14 \tr(H)K  - \tfrac14 \tr(K)H,$$
and the geodesic equation for $H(t):= g\i g_t$ becomes
$$H_t = \frac{\partial}{\partial t}|_0(g\i g_t)
= \tfrac14 \tr(HH)Id - \tfrac12 \tr(H)H.$$
 
\subheading{\nmb.{2.5} The curvature} In the setting of
\nmb!{2.1}, for vector fields $X$, $Y\in \X(\Cal N)$ and a vector
field $s:\Cal N\to T\M$ along $f:\Cal N\to \M$ we have
$$\align
R(X,Y)s &= (\nabla_{[X,Y]} - [\nabla_X,\nabla_Y])s \\
&= (K\o TK - K\o TK\o \ka_{T\M})\o T^2s\o TX\o Y,
\endalign$$
which in local coordinates reduces to the usual formula
$$R(h,k)\ell = d\Ga(h)(k,\ell) - d\Ga(k)(h,\ell) -
\Ga(h,\Ga(k,\ell)) + \Ga(k,\Ga(h,\ell)).$$
A global derivation of this formula can be found in
\cite{Kainz-Michor, 1987}
 
\proclaim{\nmb.{2.6} Proposition} The Riemannian curvature for the
canonical Riemannian metric on the manifold $\M$ of all
Riemannian metrics is given by
$$\align g\i R_g(h,k)\ell &= \tfrac14 [[H,K],L]\\
&\quad+ \frac{\dim M}{16}(\tr(KL)H-\tr(HL)K)\\
&\quad+ \frac1{16}(\tr(H)\tr(L)K-\tr(K)\tr(L)H)\\
&\quad+ \frac1{16}(\tr(K)\tr(HL)-\tr(H)\tr(KL))Id.
\endalign$$
\endproclaim
\demo{Proof} This is a long but elementary computation using
the formula from \nmb!{2.5} and
$$\align d\Ga(h)(k,\ell) &=
     -\tfrac12 k g\i h g\i \ell -\tfrac12 \ell g\i h g\i k
     -\tfrac14\tr(g\i h g\i k g\i \ell)g \\
&\quad -\tfrac14\tr(g\i k g\i h g\i \ell)g
     +\tfrac14\tr(g\i k g\i \ell)h \\
&\quad +\tfrac14\tr(g\i h g\i k)\ell
     +\tfrac14\tr(g\i h g\i \ell)k. \qed
\endalign$$
\enddemo
 
\subheading\nofrills{\bf\nmb.{2.7}. Ricci curvature }\;\;
for the Riemannian space
$(\M,G)$ {\bf does not exist}, since the mapping
$k\mapsto R_g(h,k)\ell$ is just the push forward of the section by a
certain tensor field, a differential operator of order 0.
If this is not zero, it induces a topological
linear isomorphism between certain infinite dimensional subspaces of
$T_g\M$, and is therefore never of trace class.
 
\subheading{\nmb.{2.8}. Ricci like curvature} But we may consider the
pointwise trace of the tensorial operator $k\mapsto R_g(h,k)\ell$
which we call the {\it Ricci like curvature} and denote by
$$\operatorname{Ric}_g(h,\ell)(x):= \tr(k_x\mapsto R_g(h_x,k_x)\ell_x).$$
 
\proclaim{Proposition} The Ricci like curvature of $(\M,G)$ is given by
$$\align
\operatorname{Ric}_g(h,\ell)
&= \frac{4+n(n+1)}{32}(\tr(H)\tr(L)-n\tr(HL))\\
&= -\frac n{32}(4+n(n+1))\langle h_0,\ell\rangle_ g,
\endalign$$
where $h_0:= h - \frac1n\tr(H)g$.
\endproclaim
\demo{Proof}
We compute the pointwise trace $\tr(k_x\mapsto R_g(h_x,k_x)\ell_x)$
and use the following
 
\proclaim{Lemma} For $H$, $K$, and $L\in L_{sym}(\Bbb R^n,\Bbb R^n)$
we have
$$\tr(K\mapsto [[H,K],L]) = \tr(H)\tr(L)-n\tr(HL).\qed$$
\endproclaim
\enddemo
 
We can define a $\binom02$-tensor field $\operatorname{Ric}$ on $\M$ by
$$\operatorname{Ric}(\xi,\et)(g):=
\int_M\operatorname{Ric}_g(\xi(g),\et(g))\Vol(g)$$
for $\xi$, $\et\in\X(\M)$, and $g\in\M$. Then by the proposition we
have
$$\operatorname{Ric}(\xi,\eta)=-\frac n{32}(4+n(n+1))G(\xi_0,\eta),$$
where $\xi_0$ is the vector field given by
$\xi_0(g):= \xi(g)- \frac1n\tr(g\i\xi(g))g$.
 
\subheading{\nmb.{2.9}. Scalar like curvature}
By \nmb!{2.8} there is a unique $\binom11$-tensor field
$\overline{\operatorname{Ric}}$ on $\M$ such that
$G(\overline{\operatorname{Ric}}(\xi),\eta)=\operatorname{Ric}(\xi,\eta)$
for all vector fields $\xi$, $\eta\in \X(\M)$, which is given by
$\overline{\operatorname{Ric}}(\xi)=-\frac n{32}(4+n(n+1))\xi_0$.
Again, for every  $g\in\M$ the corresponding endomorphism of $T_g\M$
is a differential operator of order 0 and is never of trace class.
But we may again form its pointwise trace as a linear vector
bundle endomorphism on $S^2T^*M\to M$, which is
a function on $M$. We call it the {\it scalar like curvature} of
$(\M,G)$ and denote it by $\operatorname{Scal}_g$.
It turns out to be the constant $c(n)$ depending only on the dimension
$n$ of $M$, because the endomorphism involved is just the projection
onto a hyperplane. We have
$$c(n)=-\frac n{32}(4+n(n+1))(\frac{n(n+1)}2-1).$$
 
\subheading{Remark} For fixed $\tilde g\in\M$ and $x\in M$ the
expression
$$G^{x,\tilde g}(h_x,k_x)
     := \langle h,k\rangle_g(x)\sqrt{\det(\tilde g\i g)(x)}$$
gives a Riemannian metric on $S^2_+T^*_xM$. It is not difficult to
see that the Ricci like curvature of $(\M,G)$ at $x$ is just the
Ricci curvature of the family of homothetic metrics on $S^2_+T^*_xM$
obtained by varying $\tilde g$, and that the scalar curvature
of $G^{x,\tilde g}$ equals the function
$c(n)/\sqrt{\det(\tilde g\i g)(x)}$ on $S^2_+T^*_xM$.
 
\heading \nmb0{3}. Analysis of the exponential mapping \endheading
 
\subheading{\nmb.{3.1}} The geodesic equation \nmb!{2.4} is an
ordinary differential equation and the evolution of
$g(t)(x)$ depends only on $g(0)(x)$ and $g_t(0)(x)$ and stays in
$S^2_+T_x^*M$ for each
$x\in M$.
 
The geodesic equation can be solved explicitly and we have
 
\proclaim{\nmb.{3.2}. Theorem} Let $g^0\in\M$ and $h\in T_{g^0}\M=\D$.
Then the geodesic in $\M$ starting at $g^0$ in the direction of
$h$ is the curve
$$g(t)=g^0 e^{(a(t)Id+b(t)H_0)},$$
where $H_0$ is the traceless part of $H:= (g^0)\i h$ (i.e.
$H_0=H-\frac{\tr(H)}n Id$) and where $a(t)$ and $b(t)\in
C^\infty(M)$ are defined as follows:
$$\align
a(t) &= \tfrac2n \log\left((1+\tfrac t4\tr(H))^2
        +\tfrac n{16}\tr(H_0^2)t^2\right)\\
b(t) &= \left\{
        {\alignedat2
        &\frac4{\sqrt{n\tr(H_0^2)}}
                \arctg\left(\frac{\sqrt{n\tr(H_0^2)}\,t}{4+t\tr(H)}\right)
                &\quad &\text{ where }\tr(H_0^2)\neq0\\
        &\frac t{1+\frac t4\tr(H)}& &\text{ where }\tr(H_0^2)=0.
        \endalignedat
        }\right.
\endalign$$
Here $\arctg$ is taken to have values in
$(-\frac\pi2,\frac\pi2)$ for the points of the manifold where
$\tr(H)\geq0$, and on a point where $\tr(H)<0$ we define
$$\arctg\left(\frac{\sqrt{n\tr(H_0^2)}\,t}{4+t\tr(H)}\right) = \left\{
        {\alignedat2
        &\arctg\text{ in }[0,\frac\pi2) &\quad &\text{ for }
                t\in [0,-\frac4{\tr(H)})\\
        &\frac\pi2 &\quad &\text{ for }t=-\frac4{\tr(H)}\\
        &\arctg\text{ in }(\frac\pi2,\pi) &\quad &\text{ for }
                t\in (-\frac4{\tr(H)},\infty).
        \endalignedat
        }\right.$$
Let $N^h:= \{x\in M: H_0(x)=0\}$, and if $N^h\neq\emptyset$ let
$t^h:=\inf\{\tr(H)(x): x\in N^h\}$. Then the geodesic $g(t)$ is
defined for $t\in [0,\infty)$ if $N^h=\emptyset$ or if $t^h\geq0$,
and it is only defined for $t\in [0,-\frac4{t^h})$ if $t^h<0$.
\endproclaim
\demo{Proof} Check that $g(t)$ is a solution of the pointwise
geodesic equation. Computations leading to this solution can be
found in \cite{Freed-Groisser, 1989}.
\qed\enddemo
 
\subheading{\nmb.{3.3}. The exponential mapping}
For $g^0\in S^2_+T_x^*M$ we consider the sets
$$\align
U_{g^0}&:= S^2T_x^*M\setminus (-\infty,-\frac4n]\,g^0,\\
L_{sym,g^0}(T_xM,T_xM) &:=
     \{\ell\in L(T_xM,T_xM): g^0(\ell X,Y)=g^0(X,\ell Y)\},\\
L^+_{sym,g^0}(T_xM,T_xM) &:=
     \{\ell\in L_{sym,g^0}(T_xM,T_xM): \ell \text{ is positive }\},\\
U'_{g^0}&:= L_{sym,g^0}(T_xM,T_xM)\setminus (-\infty,-\frac4n]\,Id_{T_xM},
\endalign$$
and the fiber bundles over $S^2_+T^*M$
$$\align
U&:= \bigcup\left\{\{g^0\}\x U_{g^0}: g^0\in S^2_+T^*M\right\},\\
L_{sym}&:= \bigcup\left\{\{g^0\}\x L_{sym,g^0}(T_xM,T_xM):
     g^0\in S^2_+T^*M\right\},\\
L^+_{sym}&:= \bigcup\left\{\{g^0\}\x L^+_{sym,g^0}(T_xM,T_xM):
     g^0\in S^2_+T^*M\right\},\\
U'&:= \bigcup\left\{\{g^0\}\x U'_{g^0}: g^0\in S^2_+T^*M\right\}.
\endalign$$
Then we consider the mapping $\Ph:U\to S^2_+T^*M$ which is given by
the following composition
$$U @>\sharp >> U'@>\ph >> L_{sym} @>\exp >> L_{sym}^+
     @>\flat >> S^2_+T^*M, $$
where $\sharp(g^0,h):= (g^0, (g^0)\i h)$ is a fiber respecting
diffeomorphism, where 
$\ph(g^0,H):= (g^0,a(1)Id + b(1)H_0)$ comes from theorem \nmb!{3.2},
where the usual exponential mapping
$\exp:L_{sym,g^0}\to L^+_{sym,g^0}$ is
a diffeomorphism (see
for instance \cite{Greub-Halperin-Vanstone, 1972, page 26})
with inverse $\log$, and where
$\flat(g^0,H):= g^0H$, a diffeomorphism for fixed $g^0$.
 
We now consider the mapping
$(pr_1,\Ph): U\to S^2_+T^*M\x_M S^2_+T^*M$.
 From the expression of $b(1)$ it is easily seen that
the image of
$(pr_1,\Ph)$ is contained in the following set:
$$\align
V:&=\left\{(g^0,g^0 \exp H): \tr(H_0^2)<\frac1n(4\pi)^2\right\}\\
&=\left\{(g^0,g):\tr\left(\left(\log((g^0)\i g)
     -\frac{\tr(\log((g^0)\i g))}n
     Id\right)^2\right)<\frac{(4\pi)^2}n\right\}
\endalign$$
Then $(pr_1,\Ph):U\to V$ is a diffeomorphism, since the
mapping $(g^0,A) \to (g^0,\ps(A))$ is an inverse to $\ph:U'\to\left\{(g^0, A):
\tr(A_0^2)<\frac1n(4\pi)^2\right\} $, where $\ps$ is given by:
  $$\ps(A) = \left\{
{\alignedat2
&\frac4n \left(e^{\frac{\tr(A)}4}
           \cos\left(\frac{\sqrt{n\tr(A_0^2)}}4\right)-1\right)Id &&\\
&\quad   +\frac4{\sqrt{n\tr(A_0^2)}}e^{\frac{\tr(A)}4}
                \sin\left(\frac{\sqrt{n\tr(A_0^2)}}4\right)A_0
&\quad&\text{ if } A_0\neq0\\
&\frac4n \left(e^{\frac{\tr(A)}4}-1\right)Id & &\text{ otherwise.}
\endalignedat
}\right.$$
 
\proclaim{\nmb.{3.4}. Theorem} In the setting of \nmb!{3.3} the
exponential mapping $\Exp_{g^0}$ is a real analytic diffeomorphism
between the open subsets
$$\gather
\Cal U_{g^0}:= \{h\in \D(S^2T^*M):(g^0,h)(M)\subset U\}\\
\Cal V_{g^0}:= \{g\in C^\infty(S^2_+T^*M):(g^0,g)(M)\subset V,
     g- g^0\in \D(S^2T^*M)\}
\endgather$$
and it is given by
$$\Exp_{g^0}(h) = \Ph\o(g^0,h).$$
The mapping $(\pi_{\M},\Exp): T\M\to \M\x \M$ is a real analytic
diffeomorphism from the open neighborhood of the zero section
$$\Cal U:= \{(g^0,h)\in C^\infty(S^2_+T^*M)\x \D(S^2T^*M):
     (g^0,h)(M)\subset U\}$$
onto the open neighborhood of the diagonal
$$\multline
\Cal V:= \{(g^0,g)\in C^\infty(S^2_+T^*M)\x C^\infty(S^2_+T^*M):
     (g^0,g)(M)\subset V, \\
     g-g^0\text{ has compact support }\}.
\endmultline$$
All these sets are maximal domains of definition for the exponential
mapping and its inverse.
\endproclaim
 
\demo{Proof}
Since $\M$ is a disjoint union of chart neighborhoods, it is
trivially a real analytic manifold, even if $M$ is not supposed to
carry a real analytic structure.
 
 From the consideration in \nmb!{3.3} it follows that $\Exp = \Ph_*$ and
$(\pi_{\M},\Exp)$ are just push forwards by smooth fiber respecting
mappings of sections of bundles. So by \cite{Michor, 1980, 8.7} they are
smooth and this applies also to their inverses.
 
To show that these mappings are real analytic, by
\cite{Kriegl-Michor, 1990} we have to check that they map real
analytic curves into real analytic curves. So we may just invoke the
description \cite{Kriegl-Michor, 1990, 7.7.2} of real analytic curves
in spaces of smooth sections:
 
{\sl For a smooth vector bundle $(E,p,M)$ a curve $c:\Bbb R\to
C^\infty(E)$ is real analytic if and only if $\hat c:\Bbb R\x M\to E$
satisfies the following condition:
\roster
\item For each $n$ there is an open neighborhood $U_n$ of
     $\Bbb R\x M$ in $\Bbb C \x M$ and a (unique) $C^n$-extension
     $\tilde c:U_n\to E_{\Bbb C}$ such that $\tilde c(\quad,x)$ is
     holomorphic for all $x\in M$.
\endroster}
In statement \therosteritem1 the space of sections $C^\infty(E)$ is
equipped with the compact $C^\infty$-topology. So we have to show
that \therosteritem1 remains true for the space $C^\infty_c(E)$ of smooth
sections with compact support with its inductive limit topology.
 
This is easily seen since we may first exchange $C^\infty(E)$ by the
closed linear subspace $C^\infty_K(E)$ of sections with support in a
fixed compact subset; then we just note that \therosteritem1 is
invariant under passing to the strict inductive limit in question.
 
Now it is clear that $\Ph$ has a fiberwise extension to a holomorphic
germ since $\Ph$ is fiber respecting from an open subset in a vector
bundle and is fiberwise a real analytic mapping. So the push forward
$\Ph_*$ maps real analytic curves to real analytic curves.
\qed\enddemo
 
\subheading{\nmb.{3.5}. Remarks}
The domain $\Cal U_{g^0}$ of definition of the exponential mapping does not
contain any ball centered at 0 for the norm derived from $G_{g^0}$.
 
Note that $\Exp_{g^0}$ is in fact defined on the set
$$\Cal U_{g^0}':= \{h\in C^\infty(S^2T^*M): (g^0,h)(M)\subset U\}$$
which is not contained in the tangent space for the
differentiable structure we use. Recall now the remarks from
\nmb!{1.2}. If we equip $\M$ with the compact
$C^\infty$-topology, then $\M$ it is not open in $C^\infty(S^2T^*M)$.
The tangent space is then the set of all tangent vectors to
curves in $\M$ which are smooth in $C^\infty(S^2T^*M)$, which is
probably not a vector space. The integral in the definition of
the canonical Riemannian metric might not converge on all these
tangent vectors.
 
The approach presented here is clean and conceptually clear, but
some of the concepts have larger domains of definition.
 
\subheading{\nmb.{3.6}. Visualizing the exponential mapping}
Let us fix a point $x\in M$ and let us consider the space
$S^2_+T_x^*M=L^+(T_xM,T_x^*M)$ of all positive definite
symmetric inner products on $T_xM$. If we fix an element
$\tilde g\in S^2_+T_x^*M$, we may define a Riemannian metric $G$ on
$S^2_+T_x^*M$ by
$$G_g(h,k) := \tr(g\i h g\i k)\sqrt{\det(\tilde g\i g)}$$
for $g\in S^2_+T_x^*M$ and $h$, $k\in T_g(S^2_+T_x^*M)=S^2T_x^*M$.
The variational method used in section \nmb!{2} leading to the
geodesic equation shows that the geodesic starting at $g^0$ in
the direction $h\in S^2T_x^*M$ is given by
$$g(t)=g^0 e^{(a(t)Id+b(t)H_0)}$$
in the setting of theorem \nmb!{3.2}.
We have the following diffeomorphisms
$$\gather
L_{sym,g^0}(T_xM,T_xM) @>\exp >> L^+_{sym,g^0}(T_xM,T_xM)
        @>{\flat=(g^0)_*}>> S^2_+T^*_xM,\\
H \mapsto \exp(H)=e^H \mapsto g^0 e^H,
\endgather$$
the manifolds $(S^2_+T_x^*M,G)$ and
$(L_{sym,g^0}(T_xM,T_xM),((g^0)_*\o \exp)^*G)$ are isometric.
For $L\in L_{sym,g^0}(T_xM,T_xM)$ we have by the affine structure
$T_LL_{sym,g^0}(T_xM,T_xM)=L_{sym,g^0}(T_xM,T_xM)$ and we get
$$((g^0)_*\o \exp)^*G)_L(H,K)
        =\sqrt{\det(\tilde g\i g^0)}\langle H,{\Cal A}_L(K)\rangle ,$$
where $\langle H,K\rangle = \tr(HK)$ and
$${\Cal A}_L =
e^{\frac12\tr(L)}
\sum_{k=0}^\infty\frac{2(\operatorname{ad}(L))^{2k}}{(2k+2)!}.$$
The geodesic on $L_{sym,g^0}(T_xM,T_xM)$ for that metric starting at 0 in the
direction $H$ is then given by $$L(t):=a(t)Id+b(t)H_0.$$
 
\midspace{12 cm}\caption{Figure 1}
 
Let us choose now a $g^0$-orthonormal basis of $T_xM$ and let
$\tilde g = n^{\frac 2n}g^0$. Then the exponential mapping at $g^0$
of $(S^2_+T^*_xM,G)$ can be viewed as the exponential mapping at 0
of $(Mat_{sym}(n,\Bbb R),\langle\quad,\quad\rangle')$, where for
symmetric matrices $A$, $B$, and $C$ we have
$$\langle A,B\rangle'_C=\frac1n \langle A,{\Cal A}_C(B)\rangle,$$
where the scale factor is chosen in such a way that we have
$\langle Id,Id\rangle'_0=1$. This exponential mapping is defined
on the set
$$\Cal U'':=\{A\in Mat_{sym}(n,\Bbb R):
        A_0\ne0\text{ or } A=\la Id \text{ with } \la>-4\}$$
by the formula
$$\multline
\exp_0(A)= \tfrac2n \log\left((1+\tfrac 14\tr(A))^2
                +\tfrac n{16}\tr(A_0^2)\right)Id \\
+\frac4{\sqrt{n\tr(A_0^2)}}
     \arctg\left(\frac{\sqrt{n\tr(A_0^2)}}{4+\tr(A)}\right)A_0.
\endmultline$$
If $A$ is traceless (i.e. $A_0=A$) and if $P$ is the plane in
$Mat_{sym}(n,\Bbb R)$ through 0, $Id$, and $A$, then
$\Exp_0(P\cap\Cal U'')\subset P$ and we can view at a
2-dimensional picture of this exponential mapping. If we
normalize $A$ in such a way that $\tr(A_0^2)=n$ the exponential
mapping is just the diffeomorphism
$$\gather
\Bbb R^2\setminus (\{0\}\x(-\infty,-\tfrac4n]) \to
        (-\frac{4\pi}n,\frac{4\pi}n)\x\Bbb R\\
\left\{{\aligned
          u(x,y) &= \tfrac4n\arctg\left(\frac{nx}{4+ny}\right)\\
          v(x,y) &= \tfrac2n \log\left(\tfrac1{16}\left((4+ny)^2
                 +n^2x^2\right)\right).
        \endaligned}\right.
\endgather$$
Here $\arctg$ is taken to have values in $(0,\pi)$ for $x\ge0$
and to have values in $(-\pi,0)$ for $x\leq0$. The images of the
straight lines and the circles can be seen in figure 1. They
correspond respectively to the images of the geodesics and to
the level sets of the distance function
$\operatorname{dist}(0,\quad)\i(r)$. For $r<\frac4n$ they are
exactly the geodesic spheres.
 
It is known that for a finite dimensional Riemannian
manifold, if for some point the exponential map is
defined in the whole tangent space, then every other
point can be joined with that by a minimizing geodesic 
(Hopf-Rinow-theorem). Here we have a nice example where,
although only a half line lacks from the domain of
definition of the exponential mapping, it is far from
being surjective.
 
\heading \nmb0{4}. Jacobi fields \endheading
 
\subheading{\nmb.{4.1}. The concept of Jacobi fields}
Let $(\Cal M,G)$ be an infinite dimensional Riemannian manifold
which admits a smooth Levi Civita connection,
and let $c:[0,a]\to\Cal M$ be a geodesic segment. By a {\it geodesic
variation} of $c$ we mean a smooth mapping
$\al:[0,a]\x(-\ep,\ep)\to\M$ such that for each fixed $s\in
(-\ep,\ep)$ the curve $t\mapsto \al(t,s)$ is a geodesic and
$\al(t,0)=c(t)$.
 
\subheading{Definition} A vector field $\xi$ along a geodesic segment
$c:[0,a]\to\Cal M$ is called a {\it Jacobi field} if and only if it is an
infinitesimal geodesic variation of $c$, i.e. if there exists a
geodesic variation $\al:[0,a]\x(-\ep,\ep)\to\M$ of $c$ such that
$\xi(t)=\dd s|_0\al(t,s)$.
 
In \nmb!{2.1} and \nmb!{2.5} we have set up all the
machinery  necessary for the usual proof that any Jacobi field
$\xi$ along a geodesic $c$ satisfies the {\it Jacobi equation}
$$\nabla_{\partial_t}\nabla_{\partial_t}\xi = R(\xi,c')c'.$$
In a finite dimensional manifold solutions of the Jacobi
equation are Jacobi fields, but in infinite dimensions one has
in general neither existence nor uniqueness of ordinary
differential equations, nor an exponential mapping.
 
Nevertheless for the manifold $(\M(M),G)$ we will show existence
and also uniqueness of solutions of the Jacobi equation for given
initial conditions, and that each solution is a Jacobi field.
 
\proclaim{\nmb.{4.2}. Lemma} Let $g(t)$ be a geodesic in
$\M(M)$. Then for a vector field $\xi$ along $g$ the Jacobi
equation has the following form:
$$\align
\xi_{tt} &=  - g_t g\i \xi  g\i g_t
     + g_t g\i \xi_t +  \xi_t g\i g_t
     + \tfrac12 \tr(g\i g_t g\i \xi_t)\,g \\
&\quad - \tfrac12 \tr(g\i g_t g\i g_t g\i \xi)\,g
     + \tfrac12 \tr(g\i g_t g\i \xi)\,g_t
     - \tfrac12 \tr(g\i \xi_t)\,g_t \\
&\quad + \tfrac14 \tr(g\i g_t g\i g_t)\,\xi
     - \tfrac12 \tr(g\i g_t)\,\xi_t.
\endalign$$
\endproclaim
 
\demo{Proof} From \nmb!{2.1} we have
$\nabla_{\partial_t}\xi = \xi_t - \Ga_g(g_t,\xi)$, thus
$$\align
\nabla_{\partial_t}\nabla_{\partial_t}\xi &= \xi_{tt} -
     \Ga(g_t,\xi)_t - \Ga(g_t,\xi_t) + \Ga(g_t,\Ga(g_t,\xi))\\
&= \xi_{tt} -d\Ga(g_t)(g_t,\xi) - \Ga(g_{tt},\xi) - 2\Ga(g_t,\xi_t)
     + \Ga(g_t,\Ga(g_t,\xi)).
\endalign$$
On the other hand we have by \nmb!{2.5}
$$R(\xi,g_t)g_t = d\Ga(\xi)(g_t,g_t) - d\Ga(g_t)(\xi,g_t) -
\Ga(\xi,\Ga(g_t,g_t)) + \Ga(g_t,\Ga(\xi,g_t)).$$
So $\xi$ satisfies the Jacobi equation if and only if
$$\xi_{tt} = d\Ga(\xi)(g_t,g_t) + 2\Ga(g_t,\xi_t).$$
By plugging in formula \nmb!{2.4} for $\Ga$ and the formula
for $d\Ga$ in the proof of \nmb!{2.6}, the result follows.
\qed\enddemo
 
\proclaim{\nmb.{4.3}. Lemma}
For any geodesic $g(t)$ and for any $k$ and $\ell\in
T_{g(0)}\M(M)$ there exists a unique vector field $\xi(t)$ along
$g(t)$ which is a solution of the Jacobi equation with
$\xi(0)=k$ and
$(\nabla_{\partial_t}\xi)(0) = \ell$.
\endproclaim
\demo{Proof} From lemma \nmb!{4.2} above we see that the Jacobi
equation is pointwise with respect to $M$, and that $\xi(t)$
satisfies the Jacobi equation if and only if at each $x\in M$ is
a Jacobi field of the associated finite dimensional manifold
treated in section \nmb!{3}. The result then follows from
the properties of the finite dimensional ordinary differential
equation involved.
\qed\enddemo
 
\proclaim{\nmb.{4.4}. Theorem}
For any geodesic $g(t)$ and for any $k$ and $\ell\in
T_{g(0)}\M(M)$ there exists a unique Jacobi field $\xi(t)$ along
$g(t)$ with initial values $\xi(0)=k$ and
$(\nabla_{\partial_t}\xi)(0) = \ell$.
 
In particular the solutions of the Jacobi equation are exactly
the Jacobi fields.
\endproclaim
\demo{Proof} We have uniqueness since Jacobi fields satisfy the
Jacobi equation and by lemma \nmb!{4.3}. Now we prove existence.
Let $h=g'(0)$ and $g^0=g(0)$. Since $k$ has compact support,
there is an $\ep>0$ such that for $s\in (-\ep,\ep)$ the tensor
field $\la(s)= g^0+sk$ is still a Riemannian metric on $M$.
Let $\tilde \ell := \ell + \Ga_{g^0}(k,h)$, and then
$W(s) := h +s \tilde\ell$ is a vector field along $\la$ which
satisfies $W(0)=h=g'(0)$ and $(\nabla_{\partial_s}W)(0) = \ell$.
Now we consider $\al(t,s):=\Exp_{\la(s)}(tW(s))$, which
is defined for all $(t,s)$ such that $(\la(t),tW(s))$ belongs to the
open set $\Cal U$ of $T\M$ defined in \nmb!{3.4}.
 
Since $h$, $k$, and $\ell$ have all compact support we have: if
$g(t)$ is defined on $[0,\infty)$, then $\ep$ can be chosen so
small that $\al(t,s)$ is defined on $[0,\infty)\x(-\ep,\ep)$;
and if $g(t)$ is defined only on $[0,-\frac4{t^h})$ then for each
$\de>0$ there is an $\ep>0$ such that $\al$ is defined on
$[0,-\frac4{t^h}-\de)\x(-\ep,\ep)$.
 
Then $J(t):= \dd s|_0\al(t,s)$ is a Jacobi field for every geodesic
segment and satisfies
$$\align
J(0) &= \dd s|_0\al(0,s) = \dd s|_0\Exp(0_{\la(s)}) = \dd s|_0\la(s) = k\\
(\nabla_{\partial_t}J)(0) &=
     \nabla_{\partial_t}|_0(t\mapsto \dd s|_0\al(t,s)) \\
&= \nabla_{\partial_s}|_0(t\mapsto \dd t|_0\al(t,s))
        =(\nabla_{\partial_s}W)(0) = \ell,
\endalign$$
where we used \nmb!{2.1}.
\qed\enddemo
 
\subheading{\nmb.{4.5}}
Since we will need it later we continue here with the
explicit expression of $\al$.
As $\la(s)\i W(s) = (Id + sK)\i(H + s\tilde L)$, where as usual
$K=(g^0)\i k$, $H=(g^0)\i h$, and $\tilde L=(g^0)\i \tilde\ell$,
from the expression in \nmb!{3.3} of the exponential mapping we get
$$\align
\al(t,s) &= \la(s) e^{Q(t,s)},\quad\text{ where }\\
Q(t,s) &= (a(t,s)Id+b(t,s)((Id + sK)\i(H + s\tilde L)-\frac{c(s)}n\,Id),
\endalign$$
where $c(s)=\tr((Id + sK)\i(H + s\tilde L))$ and where $a(t,s)$
and $b(t,s)$ are given by
$$\align
a(t,s) &= \tfrac2n\log\Bigl(\tfrac1{16}\bigl((4+tc(s))^2
            +t^2nd(s)\bigr)\Bigr),\\
b(t,s) &= \tfrac4{\sqrt{nd(s)}}
            \arctg\left(\frac{\sqrt{nd(s)}t}{4+tc(s)}\right),
                   \quad\text{where}\\
d(s) &= \tr((\la(s)\i W(s))^2_0)=f(s)-\frac{c(s)^2}n,\quad\text{and} \\
f(s) &= \tr\left(((Id + sK)\i(H + s\tilde L))^2\right).
\endalign$$
In fact $b(t,s)$ should be defined with the same care as in the
explicit formula for the geodesics in \nmb!{3.2}. This is
omitted here.
 
Before computing the Jacobi fields let us introduce some notation.
For every point in $M$ the mapping $(H,K)\mapsto \tr(HK)$
is an inner product, thus the quadrilinear mapping
$$T(H,K,L,N):= \tr(HL)\tr(KN)-\tr(HN)\tr(KL)$$
is an algebraic curvature tensor, \cite{Kobayashi-Nomizu, I,
page 198}. We will also use
$$S(H,K):= T(H,K,H,K) = \tr(H^2)\tr(K^2)-\tr(HK)^2.$$
 
Let $P(t):= g(t)\i g'(t)$ then it is easy to see from \nmb!{3.2} that
$$P(t)= e^{-\frac12 na(t)}
     \left(\frac{4\tr(H)+nt\tr(H^2)}{4n}Id + H_0\right).$$
We denote by $P(t)^\bot$ the $\binom11$-tensor field of trace
$e^{-\frac12 na(t)}$ in the plane through $0$, $Id$, and $H_0$ which
at each point is orthogonal to $P(t)$ with respect to the inner
product $\tr(HK)$ from above.
 
\proclaim{\nmb.{4.6}. Lemma} In the setting of \nmb!{4.5} we have
$$\multline
\dd s|_0 Q(t,s) = \frac{\tr(H\hat L)}{\tr(H^2)}tP(t)
+ \frac{T(\hat L,H,Id,H)}{\tr(H^2)}tP(t)^\bot \\
+ b(t)\left(-\frac{T(H,Id,\hat L,Id)}{S(H,Id)}H_0 + \hat L_0\right),
\endmultline$$
where $\hat L := -KH + \tilde L$.
\endproclaim
 
\demo{Proof}
 From the expression of $Q(t,s)$ we have
$$\dd s|_0 Q(t,s) = \left(\dd s|_0a(t,s)Id
        + \dd s|_0b(t,s)H_0 + b(t,0)\hat L_0\right).$$
Now,
$$\align
\dd s|_0a(t,s) &= \frac2n\frac{8tc'(0)+nt^2f'(0)}{16+8tc(0)+nt^2f(0)}\\
\dd s|_0b(t,s) &= -\frac{d'(0)}{2d(0)}\,b(t)
     +\frac2{d(0)}\frac{(4+tc(0))td'(0)-2t^2c'(0)d(0)}{16+8tc(0)+nt^2f(0)}.
\endalign$$
 From the definitions of $c$, $d$, and $f$ we have
$$\align
c(0) &= \tr(H), \qquad c'(0) = \tr(\hat L) \\
d(0) &= \tr(H_0^2)  =\tfrac1nS(H,Id), \\
f(0) &= \tr(H^2), \qquad f'(0) = 2\tr(H\hat L), \\
d'(0) & = \tfrac2n\,T(H,Id,\hat L,Id)
\endalign$$
and $16+8tc(0)+nt^2f(0) = 16\,e^{\frac{na(t)}2}$, and so we get
$$\align
\dd s|_0a(t,s) &=
\frac1{4n}e^{-\frac12na(t)}(4t\tr(\hat L)+nt^2\tr(H\hat L)), \\
\dd s|_0b(t,s) &= -\frac{T(H,Id,\hat L,Id)}{S(H,Id)}b(t) \\
+ &\frac{e^{-\frac12na(t)}t}{4S(H,Id)}
     \bigl((4 + t\tr(H))T(H,Id,\hat L,Id)-t\tr(\hat L)S(H,Id)\bigr).
\endalign$$
If we collect all terms and compute a while we get the result.
\qed\enddemo
 
\proclaim{\nmb.{4.7}. Theorem} Let $g(t)$ be the geodesic in
$(\M(M),G)$ starting from $g^0$ in the direction $h\in T_{g^0}\M(M)$.
For each $k,\ell\in T_{g^0}\M(M)$
the Jacobi field $J(t)$ along $g(t)$ with initial conditions
$J(0)=k$ and $(\nabla_{\partial_t}J)(0) = \ell$ is given by
$$\align J(t) &=
\frac{\tr(H\hat L)}{\tr(H^2)}\,t\,g'(t) + \frac{\tr(\hat L)\tr(H^2)
     - \tr(\hat LH)\tr(H)}{\tr(H^2)}\,t\,(g'(t))^\bot \\
&\quad+ b(t)g(t)\left(-\frac{\tr(H_0\hat L_0)}{\tr(H_0^2)}H_0 + 
     \hat L_0\right)\\ 
&\quad+ g(t)\sum_{m=1}^\infty\frac{(-\operatorname{ad}(b(t)H))^m}{(m+1)!}
(b(t)\hat L) + k(g^0)\i g(t),
\endalign$$
where $H=(g^0)\i h$, $K=(g^0)\i k$, $L=(g^0)\i \ell$,\newline
$\hat L=-KH+L+(g^0)\i\Ga_{g^0}(k,h)$,
and $(g'(t))^\bot = g(t)P(t)^\bot $.
\endproclaim
 
\demo{Proof}
By theorem \nmb!{4.4} we have
$$\align
J(t) &= \dd s|_0\al(t,s)\qquad\text{ and then }\\
&= g^0 \dd s|_0 e^{Q(t,s)} + k\,e^{Q(t,0)}\\
&= g^0 e^{Q(t,0)}\left(
     \sum_{m=0}^\infty\frac{(-\operatorname{ad}(b(t)H))^m}{(m+1)!}
     (\dd s|_0 Q(t,s))\right) + k\,e^{Q(t,0)}.
\endalign$$
The result now follows from lemma \nmb!{4.6}
\qed\enddemo

\enddocument